%% file: main.tex
\documentclass[letterpaper, 10 pt, conference]{ieeeconf}  % Comment this line out if you need a4paper

\IEEEoverridecommandlockouts                              % This command is only needed if 

\title{\LARGE 
Inner Approximations of the Positive-Semidefinite Cone via Grassmannian Packings}
\author{Tianqi Zheng, James Guthrie, and Enrique Mallada\thanks{T. Zheng, J. Guthrie, and E. Mallada are with the Department of Electrical and Computer Engineering, Johns Hopkins University, Baltimore, MD 21218, USA. Email: {\tt\small \{tzheng8, jguthri6 ,mallada\}@jhu.edu}. 
}}

\usepackage{amsthm}
\usepackage{amsmath}
\usepackage{amsfonts}
\usepackage{amssymb}
\usepackage{mathtools}

\allowdisplaybreaks
\usepackage{subfiles}
\usepackage{environ}
\usepackage{accents}
\usepackage{balance}
\usepackage{lineno}
\usepackage{graphicx}
\usepackage{subfig}
\usepackage{cite}
\newtheorem{thm}{Theorem}

\newtheorem{prop}[thm]{Proposition}

\newtheorem{defn}{Definition}

\newtheorem{prob}{Problem}

\allowdisplaybreaks

\let\labelindent\relax 
\usepackage{enumitem}
\setlist[enumerate]{leftmargin=*}
\setlist[itemize]{leftmargin=*}
 
\input{edits.tex}

\DeclareSymbolFont{bbold}{U}{bbold}{m}{n}
\DeclareSymbolFontAlphabet{\mathbbold}{bbold}

\newcommand{\eucd}{\mathbb{R}}

\begin{document}
\maketitle
\begin{abstract}
    We investigate the problem of finding inner approximations of positive semidefinite (PSD) cones. We develop a novel decomposition framework of the PSD cone by means of conical combinations of smaller dimensional sub-cones. We show that many inner approximation techniques could be summarized within this framework, including the set of (scaled) diagonally dominant matrices, Factor-width $k$ matrices, and Chordal Sparse matrices. Furthermore, we provide a more flexible family of inner approximations of the PSD cone, where we aim to arrange the sub-cones so that they are maximally separated from each other. In doing so, these approximations tend to occupy large fractions of the volume of the PSD cone. The proposed approach is connected to a classical packing problem in Riemannian Geometry. Precisely, we show that the problem of finding maximally distant sub-cones in an ambient PSD cone is equivalent to the problem of packing sub-spaces in a  Grassmannian Manifold. We further leverage existing computational method for constructing packings in Grassmannian manifolds to build tighter approximations of the PSD cone. Numerical experiments show how the proposed framework can balance between accuracy and computational complexity, to efficiently solve positive-semidefinite programs.
\end{abstract}


\section{Introduction}

Semidefinite programs (SDPs) are a class of convex optimization problems whose domain is formed by the intersection of the cone of positive semidefinite (PSD) matrices with affine subspaces. Having high expressive power, SDPs have been studied extensively over the past few decades \cite{blekherman2012semidefinite,vandenberghe1996semidefinite}, and have enable many applications in  robust control \cite{parrilo2000structured,8431617}, trajectory planning \cite{deits2015computing} and robotics \cite{5717711}.  In theory, SDPs could be solved to arbitrary accuracy in polynomial time via efficient algorithms, such as  second order interior point methods (IPMs). However, in practice, when the PSD cone, i.e., the set of all PSD matrices, is high dimensional, the solving time and memory requirements tend to be impractical \cite{benson2000solving}. Over the past few years, extensive research has been performed with the aim of to improving the scalability of semidefinite programs \cite{majumdar2019survey}. Since the computational cost of solving an SDP is closely related to the size of the largest semidefinite constraint, several methods have been proposed to decompose, or approximate high-dimensional PSD constraints into conic combinations of smaller PSD cones~\cite{ahmadi2019dsos,zheng2019block,blekherman2020sparse, song2021approximations}.

This has led to several methods that provide inner approximations for the PSD cone, some of which may still lead to the same SDP solution. Examples include, the set of \textit{diagonally dominant} ($DD_n$) matrices and  \textit{scaled diagonally dominant} ($SDD_n$) matrices, which can be written as conic combinations of rank one PSD matrices and conic combinations of second order cones, respectively~\cite{ahmadi2019dsos}, and lead to approximations that can be solved using
% resp. sdd) matrices can be equivalently written as a number of linear constraints (resp. second order cone constraints . 
% In this way, SDPs get reduced to optimization problems over the cone of dd matrices (resp. sdd matrices) which can be solved using 
linear programs (LP) or resp. second-order cone programs (SOCP). 
Along the same line, \textit{Factor-width $k$} matrices, which can be viewed as a generalization of scaled diagonally dominant matrices~\cite{boman2005factor}, can be decomposed as a combination of PSD matrices of rank at most $k$. Finally, one could exploit chordal sparsity patterns of matrices \cite{vandenberghe2015chordal}, i.e., sparsity induced by a chordal graph, to equivalently express a large semidefinite constraint as a set of smaller semidefinite constraints. Unfortunately, the scalability of chordal decompositions is limited by the size of the largest clique.
Moreover,  approximating the PSD cone with $ DD_n$ or $ SDD_n$ is conservative, which means that either feasible solutions of the constrained problem may be sub-optimal or that no solution exists. Further, enforcing these constraints is problematic as $n$ increases, since the number of constraints exhibit a combinatorial growth~\cite{ahmadi2019dsos, zheng2019block}.

Our work mostly aligns with the efforts on obtaining inner approximations of the PSD cone, but proposes a more flexible family of inner approximations based on Grassmannian packing problems. More specifically, we develop a novel decomposition framework of the PSD cone that generalizes the above-mentioned  approximation techniques. Using this decomposition, we develop novel inner approximations of the PSD cone, based on conic combinations of smaller dimensional sub-cones, such that they are maximally separated from each other. In doing so, these approximations tend to occupy large fractions of the volume of the PSD cone and provide more accurate inner approximations. 
Such problem of placing maximally distant sub-cones is closely related to the classical packing problem in Grassmannian manifolds. The Grassmannian packing problem~\cite{conway1996packing} aims to answer the following question: how should $N$ $k$-dimensional subspaces of an $n$-dimensional Euclidean space be arranged so that they are as far apart as possible? The problem of finding Grassmannian packings is a fundamental problem of geometry and has applications in many fields such as signal processing \cite{kutyniok2009robust} , information theory \cite{love2003grassmannian,gohary2009noncoherent}, and wireless communications \cite{strohmer2003grassmannian}. The hope is that by establishing these connections, this work will ignite a new search for better approximations that exploit classical tools from, e.g., information theory.

The rest of the paper is organized as follows. Section \ref{sec:preliminary} introduces some preliminary terminology on packing problems in Grassmannian manifolds. Section \ref{sec: General framework Literature} reviews several existing inner approximations of the PSD cone and provide a novel decomposition framework of the PSD cone which unifies existing examples. Section \ref{sec:Volume} further provides our Grassmannian packing based approximations, by connecting distances among elements in the Grassmannian manifold with a properly defined distance between sub-cones.
% develop novel inner approximations of the PSD cone, based on conical combinations of smaller dimensional sub-cones, where we aim to pack these sub-cones so that they are as far apart possible, on a properly defined metric. In doing so, we introduce a notion of cone distance between two sub-cones and showed that the problem of packing sub-cones in an ambient PSD cone is equivalent to a packing problem in Grassmannian manifolds. 
Numerical examples, that leverage existing computational methods for packings in Grassmannian~\cite{dhillon2008constructing}, illustrate how one can balance between accuracy and computational complexity using the proposed approach. 
% The first one compares the inner approximation accuracy of scaled diagonally dominant (sdd) matrices, Factor-width 3 matrices ($\mathcal{FW}_n^3$) with our approach. The second numerical example studies how many sub-cones are needed to approximate the ambient PSD cone with high accuracy.
We finally conclude in Section \ref{sec:conlusions}.

\textit{Notation}: Denote the set of $ n \times n$ real symmetric matrices by $\mathbb{S}^n$. A matrix $A \in \mathbb{S}^n$ is positive semi-definite (PSD), denoted by $A \succeq 0 $, if $ x^T Ax \geq 0, \forall x \in \eucd^n$. We use standard notation $\mathbb{S}^n_+ $ to denote the set of PSD matrices and $ \langle \cdot , \cdot \rangle$ denotes the inner product in the appropriate space. The set of $n \times n$ orthogonal matrices is denoted by $O(n) $. Given a set $ S$, its cardinality is denoted by $ |S|$. The set $\{1, \dots, n \}$ is denoted by $ [n]$. The Frobenius norm of a Matrix $X$ is denoted and calculated as $\|X \|_F^2 = \mathrm{ trace}(X^T X) $, where the trace operator sums the diagonal entries of the matrix. The conic hull for a set $C$, i.e., the set of all conic combinations of points in  $C$, is given by  $  \mathrm{cone}(C) := \{\sum_{i=1}^k \alpha_ix_i: x_i \in C, \; \alpha_i \in \eucd_{\geq 0}, \; k \in \mathbb{N} \}$. Given a matrix $ F \in \eucd^{n \times n}$ and a non-empty index set $S \subseteq [n]$,  $F_S \in \eucd^{n \times |S|}$ denotes the sub-matrix of $F$ comprised by columns indexed by $S$.

\section{Preliminary} \label{sec:preliminary}

\subsection{Packing in Grassmannian Manifolds}
In this section we formalize the problem of finding maximally distant $k$-dimensional subspaces of real Euclidean $n-$dimensional space $ \eucd^n$.
We thus consider Grassmannian manifold  $G(k, n) $, i.e., the set of all $k-$dimensional subspaces of $\eucd^n$. 
% Over the past several decades, many researchers have investigated the following packing problem in Grassmannian manifolds: how should $N$ $k-$dimensional subspaces of $n-$dimensional Euclidean space be arranged so that they are as far apart possible \cite{conway1996packing} ? 
% In order to formally define this Grassmannian packing problem, we first introduce a useful metric in Grassmannian manifolds: the chordal distance. 
In order to define a proper notion of distance between elements in $G(k,n)$ it is useful to consider the principal angles between subspaces.
\begin{defn}[Principal Angles \cite{conway1996packing}]
Suppose that  $ \mathcal{S},\mathcal{T} $ are two $k-$subspaces in $ G(k, n)$. The principal angles $ \theta _1, \dots, \theta_k \in [0, \pi/2]$ between  $ \mathcal{S} $ and $ \mathcal{T} $ are defined by:
\begin{align*}
    \cos \theta_i = \max_{ u \in  \mathcal{S}} \max_{ v \in  \mathcal{T}} u \cdot v = u_i \cdot v_i,
\end{align*}
for $i \in [k]$, subject to $ u \cdot u =v \cdot v = 1,u \cdot u_j = 0, v \cdot v_j = 0 \;(1 \leq j \leq i-1)  $.
\end{defn}

A more computational definition of the principal angles based on  singular value decomposition (SVD) is provided in \cite{bjorck1973numerical}, which we describe next. Suppose that columns of $ S$ and $T$ form orthonormal bases for $ \mathcal{S} , \mathcal{T} \in G(k, n) $ respectively. Therefore, we have $S \in \eucd^{n \times k} $ such that $ S^T S = I_k , \mathrm{range} (S) = \mathcal{S} $, and analogously for $T$.\footnote{Note that this is equivalent to say that $ S$ belongs to the Stiefel manifold $ V_k(\eucd^n):= \{A \in \eucd^{n \times k} : A^T A= I_k\}$.} 
% , i.e., the set of ordered orthonormal $k$-tuples of vectors in $ \eucd^n$.
Next, we compute the singular value decomposition of the product $ S^T T = U  \Sigma V^T$, where {$ U, V $ are $ k \times k$ unitary} matrices and $\Sigma $ is a nonnegative, diagonal matrix with non-increasing entries. The matrix $\Sigma$ of singular values is uniquely determined, with entries being cosines of the principal angles between  $ \mathcal{S} $ and $ \mathcal{T}$: $\Sigma_{ii}   = \cos \theta _i , i \in [k]$.
Grassmannian manifolds admit many metrics  such as the spectral distance, the Fubini-Study distance, the geodesic distance etc. In this paper, we focus on the chordal distance, which is the easiest to work with and has a number of desirable features such as its square is differentiable everywhere.
\begin{defn}[Chordal Distance \cite{conway1996packing}]
Suppose we have  two $k-$dimensional subspaces $ \mathcal{S} , \mathcal{T} $ in $ G(k, n) $ and the columns of $ S ,T \in \eucd^{n \times k}$ form orthonormal bases for $ \mathcal{S} , \mathcal{T}$, respectively.
The chordal distance between $ \mathcal{S} $ and $ \mathcal{T}$, as a function of $(S ,T)$, is given by
\begin{align} \label{eq:Definition_Chordal_Distance}
    d_{\mathrm{ chord}}(S , T ) :=& \sqrt{ \sin^2\theta_1 + \dots + \sin^2\theta_k}  \nonumber \\
    =& [k - \| S^T T \|_F^2 ]^{1/2},
\end{align}
where  $ \theta _1, \dots, \theta_k \in [0, \pi/2]$ denotes the principal angles  between  $ \mathcal{S} $ and $ \mathcal{T} $ .
\end{defn}
Note that based on the SVD decomposition, we have $ S^T T = U  \Sigma V^T$, where $\Sigma_{ii}   = \cos \theta _i , i \in [k]$, and the second equation comes from the fact that
$  \| S^T T \|_F^2 = \sum_{i=1}^k \cos^2\theta_k  $. Now, we formally define the packing problem in Grassmannian manifolds with chordal distance.

\begin{prob}[Packing problem in Grassmannian manifolds with chordal distance] \label{Prob:Grass_Packing}
Given the Grassmannian manifold $G(k, n) $ of $k-$dimensional subspaces of the real Euclidean $n-$dimensional space $ \eucd^n$, find a set of N $k-$dimensional subspaces $ \{\mathcal{S}_1,\dots,\mathcal{S}_N \} \subseteq G(k, n) $  spanned by the matrices $\mathcal{F}=\{F_1,\dots,F_N \}$, that solves the mathematical program
\begin{align}
    \max_{\mathcal{F}: | \mathcal{F}|=N } \;\; \min_{F_i,F_j \in \mathcal{F}, i\neq j } d_{\mathrm{chord}}(F_i,F_j).
\end{align}
\end{prob}

\section{PSD Cone Decomposition Framework}\label{sec: General framework Literature}
Given $b \in \eucd^m$ and matrices $C, A_1,\dots,A_m \in \mathbb{S}^n   $, the standard primal form of a semidefinite program (SDP) is
\begin{align}\label{eq:SDP_form}
    \min_{X}\;\; & \langle C,X\rangle  \nonumber \\
    s.t.\;\; &  \langle A_i,X\rangle = b_i\;\; , i \in [m] \\
    &X \in \mathbb{S}^n_+, \nonumber
\end{align}
where the optimization variable $X$ is constrained to be positive semidefinite (psd). That is, linear optimization over the cone of PSD matrices can be addressed by a semidefinite program. Although many algorithms, such as interior point methods, could solve SDPs to arbitrary accuracy in polynomial time, scaling as $O(n^2m^2+n^3m) $ per iteration \cite{alizadeh1995interior}. In practice,  when the PSD cone is high dimensional, the solving time and memory requirements tend to be impractical. More specifically, when $m$ is fixed, the computational cost of solving an SDP is primarily determined by the size of the largest semidefinite constraint $n$. 
This motivates the study of methods to decompose, or approximate high-dimensional PSD constraints into conic
combinations of smaller PSD cones. Such alternatives allow one to trade off computation time with solution quality.
In this section, we first review several inner approximation techniques of the PSD cone $\mathbb{S}^n_+ $. Then we provide a novel PSD cone decomposition framework by means of conical combinations of smaller dimensional subcones to unify these results.

\subsection{Inner Approximations of PSD cone}
The first example is the set of diagonally dominant matrices and scaled diagonally dominant matrices. 

\begin{defn}[$DD_n$ and $SDD_n$ matrices]
A symmetric matrix $A=[a_{ij}]$  is diagonally dominant (dd) if $a_{ii} \geq \sum _{j \neq i} |a_{ij}|$. A symmetric matrix A is scaled diagonally dominant (sdd) if there exists a diagonal matrix $D$, with positive diagonal entries, such that $DAD$ is dd. We denote the set of $n \times n$ dd (resp. sdd) matrices with $DD_n \;(\text{resp.}~SDD_n)$. Note that $ DD_n \subseteq SDD_n \subseteq \mathbb{S}^n_+$.
\end{defn}
According to the extreme ray characterization of diagonally dominant matrices by Barker and Carlson \cite{barker1975cones}, a symmetric matrix $M \in DD_n$ if and only if it can be written as $ M = \sum_{i=1}^{n^2} \alpha_i v_iv_i^T, \; \alpha_i \in \mathbb{S}^1_+$, where $ \{ v_i\}$ is the set of all nonzero vectors in $\eucd^n$ with at most 2 nonzero components, each equal to $ \pm 1$. The vectors $ \{ v_i\}$ are the \textit{extreme rays} of the $ DD_n $ cone.  A similar decomposition is available for scaled diagonally dominant matrices. According to  \cite{boman2005factor}, we know that any scaled diagonally dominant matrix $M$ can be written as $ M = \sum _{i=1}^{\binom{n}{2}} V_i \Lambda_i V_i^T $, where $V_i$ is an $n \times 2 $ matrix whose columns each contain one nonzero element which is equal to 1, and $\Lambda_i \in \mathbb{S}^2_+$.

The cone of $ DD_n$ (resp. $SDD_n$) can be equivalently written as a number of linear constraints (resp. second order cone constraints), and thus linear optimization over the set of diagonally dominant matrices (resp.scaled diagonally dominant matrices) can be addressed using linear programming (LP) (resp. second order cone programming (SOCP) ) \cite{ahmadi2019dsos}. Working with these classes of convex optimization problems allows one to take advantage of high-performance LP and SOCP solvers. However, as already pointed in \cite{ahmadi2019dsos}, approximating the PSD cone by the set of scaled diagonally dominant matrices $ SDD_n$ is conservative, which means the restricted problem may be infeasible or the optimal solution may be significantly different from that of the original SDP. Furthermore, enforcing these constraints could be problematic as $n$ increases, since the number of constraints grows in a combinatorial fashion  $\binom{n}{2}$.

A relevant generalization of the set of scaled diagonally dominant matrices is the set of factor-width k matrices. 

\begin{defn}[Factor-width k matrices \cite{boman2005factor}]
The factor width of a symmetric matrix $A$ is the smallest integer $k$ such that there exists a matrix $V$ where $ A = VV^T$ and each column of $ V$ contains at most $k$ non-zeros. We denote the cone of $n \times n$ symmetric matrices of factor width no greater than $k$ by $ \mathcal{FW}_n^k$.
\end{defn}
By definition, we have  that $ \mathcal{FW}_n^k \subseteq \mathbb{S}^n_+ $ for all $k \in [n]$ and $\mathcal{FW}_n^2 = SDD_n$, and $ \mathcal{FW}_n^n = \mathbb{S}^n_+ $ \cite{ahmadi2019dsos}. $ \mathcal{FW}_n^k$ can be decomposed into a sum of PSD matrices of rank at most $k$ as follows: $Z \in \mathcal{FW}_n^k$ if and only if it can be written as $ Z = \sum _{i=1}^{\binom{n}{k}} V_i \Lambda_i V_i^T $, where $V_i$ is an $n \times k $ matrix whose columns each contain one nonzero element which is equal to 1, and $\Lambda_i \in \mathbb{S}^k_+$~\cite{zheng2019block}. One could transform an optimization problem over the PSD cone into an conic optimization problem of smaller dimension. However, enforcing this constraint is problematic as $n$ or $k$ increases, which grows in a combinatorial fashion  ${\binom{n}{k}} $. Therefore, using factor-width k matrices to approximate the PSD cone is not trivial in a practical way.

Lastly, it is often important to exploit sparsity in large semidefinite optimization problem and chordal graph properties from graph and sparse matrix theory plays a central role. The following proposition is the primary result that exploit the chordal sparsity and decompose a larger semidefinite constraint as a set of smaller semidefinite constraints. Details could be found in the survey of chordal graphs and semidefinite optimization \cite{vandenberghe2015chordal}. 

\begin{prop}[Chordal sparse matrices \cite{vandenberghe2015chordal}]\label{Prop:Chordal} 
 Let $ \mathcal{G}$ be an undirected chordal graph with edges $\mathcal{E}$, vertices $\mathcal{V}$ and a set of maximal cliques $\mathcal{T} = \{C_1,\dots,C_p \}$. The set of matrices with a sparsity pattern defined by $ \mathcal{G}$ is defined as: $\mathbb{S}^n(\mathcal{E} ) := \{ X \in \mathbb{S}^n : X_{ij} = 0 \;\mathrm{if}\; (i,j) \not\in \mathcal{E} \; \mathrm{for}\; i \neq j \} $. Then $X \in \mathbb{S}^n(\mathcal{E} ) $ is positive semidefinite if and only if $X = \sum^p_{k=1} X_k ,\; X_k \succeq 0, \; k \in [p]$, where $X_k \in S( C_k ) := \{ X \in \mathbb{S}^n : X_{ij} = 0 \;\mathrm{if}\; (i,j) \not\in C_k \times C_k \} $.
\end{prop}

\subsection{Decomposition of the PSD Cone}
In the above section, we list three existing techniques of inner approximation for the PSD cone, by means of decomposing a larger PSD matrices by conic combinations of smaller dimensional sub-cones, with additional equality or inequality constraints. Motivated by this idea, we propose a simple, yet insightful, decomposition  of the PSD cone, with the hope to render new methods for approximating $\mathbb{S}^n_+$.

\begin{thm}  \label{lem:General PSD Decomposition Framework}
Consider a collection $ \mathcal{S}$ of subsets of $[n]$, i.e., if an index set $S \in \mathcal{S}$, then $ S \subseteq  [n]$. Then we can decompose the positive semidefinite cone $\mathbb{S}^n_+ $ as the conic hull of
the (infinite) union of lower dimensional cones, that is,
\begin{align} \label{eq:GeneralFramwwork_Decomposition_of_PSD_Cone_}
    \mathbb{S}^n_+ = \mathrm{cone}( \bigcup_{S \in \mathcal{S} } \bigcup_{ F \in O(n)} F_S[\mathbb{S}^{|S| }_+ ]  F_S^T  ),
\end{align}
where the set
\begin{align*}
    F_S[\mathbb{S}^{|S| }_+ ]  F_S^T := \{ X \in \mathbb{S}^n_+ :  X =  F_S Y  F_S^T, Y \in \mathbb{S}^{|S| }_+  \}
\end{align*}
is in and of itself a sub-cone. 

\textit{Proof}: See Appendix~\ref{appendix:A}
\end{thm}

As a result of Theorem \ref{lem:General PSD Decomposition Framework}, the PSD cone could be viewed as conic hull of sub-cones, possibly of different dimensions. Therefore, an inner approximation of PSD cone following the structure of decomposition \eqref{eq:GeneralFramwwork_Decomposition_of_PSD_Cone_} is given by: choosing specific choices of $\mathcal{S} $ and a finite subset $ \mathcal{F} \subseteq O(n) $ ,i.e., PSD matrices $X$ of the form
\begin{align}
    X = \sum_{S \in \mathcal{S}} \sum _{F \in \mathcal{F}}F_S X_{F,S}F_S^T,  \;\; X_{F,S} \in \mathbb{S}^{|S| }_+ \label{eq:novel_InnerApproximation_PSDCone}
\end{align}
Next, we illustrate how the set of (scaled) diagonally dominant matrices, factor-width $k$ matrices and chordal sparse matrices could be construed as in \eqref{eq:novel_InnerApproximation_PSDCone}.

\begin{enumerate}
    \item The set of diagonally dominant matrices, $DD_n$. As mentioned in above subsection, a symmetric matrix $M \in DD_n$ if and only if it can be written as $ M = \sum_{i=1}^{n^2} \alpha_i v_iv_i^T, \; \alpha_i \in \mathbb{S}^1$, where $ \{ v_i\}$ is the set of all nonzero vectors in $\eucd^n$ with at most 2 nonzero components, each equal to $ \pm 1$.
    We could choose $ \mathcal{S} = \{S \subseteq [n]: | S| \leq 2 \}$ and $ \mathcal{F} = \{I_n^{\pm}\} \subseteq O(n)$, where $I_n^{\pm}$ being the set of all diagonal matrices with elements $\pm 1$ in its diagonal. $| S| = 2 $ corresponds to the case where  $v_i$ has 2 nonzero components.
    
    \item $SDD_n(\mathcal{FW}_n^2)$ and $\mathcal{FW}_n^k$. Recall that $Z \in \mathcal{FW}_n^k$ if and only if it can be written as $ Z = \sum _{i=1}^{\binom{n}{2}} V_i \Lambda_i V_i^T $, where $V_i$ is an $n \times k $ matrix whose columns each contain one nonzero element which is equal to 1, and $\Lambda_i \in \mathbb{S}^k$.
    Consequently, for the set of matrices $\mathcal{FW}_n^k$, one could choose $ \mathcal{S} = \{S \subseteq [n]: | S| \leq k \}$ and $ \mathcal{F} = \{I_n\}$.
    \item Chordal Sparse Matrices. According to Proposition \ref{Prop:Chordal}, one could choose the set $ \mathcal{S} = \{C_k\}, k\in [p]$, i.e., the collection of index set of all the maximal cliques of the graph $\mathcal{G}(\mathcal{V},\mathcal{E})$, and $ \mathcal{F} = \{I_n\}$.
\end{enumerate}

\section{Maximally Distant Frames and Grassmannian Packings}\label{sec:Volume}
The examples in the previous section showed that existing approximation techniques could be considered as special instance of the form \eqref{eq:novel_InnerApproximation_PSDCone}. However, these approaches restrict a priori the relation between the set $  \mathcal{F}$ and $  \mathcal{S}$ (also their cardinalities). In this section, we provide a more flexible family of inner approximations by focusing on Grassmannian packings.

In the rest of the paper, we restrict our attention to the case where $S $ is a fixed subset of cardinality $s$, i.e., the sub-cones are of the same size, and the matrices $F \in \eucd^{n \times s}$ are tall matrices composed by $s $ orthonormal vectors, i.e., $F$ belongs to the Stiefel manifold $V_s(\eucd^n )$. In doing so, we will further simplify the decomposition  \eqref{eq:GeneralFramwwork_Decomposition_of_PSD_Cone_} and provide novel inner approximations of the PSD cone, which aimed to arrange the sub-cones such they are maximally separated from each other and tend to occupy large fractions of the volume of ambient PSD cone. The proposed approach is rooted in the classical packing problem in Grassmannian manifolds.

\begin{thm} \label{lem:simplified_Decomposition}
Consider the case $ S \subseteq  [n]$ is a fixed subset of cardinality $s$ and the matrices $F \in \eucd^{n \times s}$ are tall matrices composed by $s $ orthonormal vectors, i.e., $ F \in V_s(\eucd^n )$ . Then, we could simplify the cone decomposition \eqref{eq:GeneralFramwwork_Decomposition_of_PSD_Cone_} to
\begin{align}\label{eq:simplified_PSD_Decomposition}
    \mathbb{S}^n_+ = \mathrm{cone}(  \bigcup_{ F \in V_s(\eucd^n )} F[\mathbb{S}^s_+ ]  F^T  ) 
\end{align}
where the set 
\begin{align}\label{eq:sub-cones}
    F[\mathbb{S}^s_+ ]  F^T := \{ X \in \mathbb{S}^n_+ :  X =  F Y  F^T, Y \in \mathbb{S}^{s}_+  \}
\end{align}
is in and of itself a sub-cone. 

\textit{Proof}: See Appendix~\ref{appendix:B}
\end{thm}

The above Theorem and the simplified decomposition \eqref{eq:simplified_PSD_Decomposition} illustrate an important fact: Since one could fully recover the PSD cone $\mathbb{S}^n_+$ by conical combinations of elements drawn from distinct sub-cones of the form $F[\mathbb{S}^s_+ ]  F^T  $ spanned by a set of frames $F \in V_s(\eucd^n)$. Therefore, when selecting a finite set of frames $ \mathcal{F} = \{ F_k\} \subseteq V_s(\eucd^n)$ to generate an inner approximation of PSD cone, one should seek to select the frames $ \{ F_k\}$ such that the corresponding sub-cones $C_s(F_k):= F_k[ \mathbb{S}_+^s]F_k^T$ are placed as far from each other as possible with in $\mathbb{S}^n_+ $, such that these approximations tend to occupy large fractions of the volume of the PSD cone. 
Precisely, we have the following form of inner approximation of the PSD cone, 
\begin{align} \label{eq:Final_approximation}
   \sum_{k=1}^N F_k[ \mathbb{S}_+^s]F_k^T
\end{align}
for a fixed number of N frames $ \mathcal{F} = \{ F_k\}, k\in [N] $. We aim to construct these frames such that  the corresponding sub-cones $C_s(F_k)$ are maximally distant from each other.

Consequently, we seek to solve the following problem: \emph{ How should N $s-$dimensional sub-cones $C_s(F_k), \; k\in [N]$, of $n-$dimensional PSD cone $\mathbb{S}^n_+$ be arranged so that they are as far apart possible?}

To make sense of this statement, we first define a notion of cone separation between two sub-cones $C_s(F_i) =F_i[ \mathbb{S}_+^s]F_i^T $ and $C_s(F_j)=F_j[ \mathbb{S}_+^s]F_j^T$. 
Though there are many alternative notions of separation between sub-cones that could employed, we focus here on the distance between the sub-cone central rays given by $F_i[\alpha_i I_s]F_i^T$, $\alpha_i\geq0$. Further, since all sub-cones share the origin, normalization (setting $\alpha_i=1$) is necessary. This leads to the following notion of cone distance.

\begin{defn}[Cone distance]
Consider two Frames $F_i, F_j \in V_s(\eucd^n )$ and the corresponding sub-cones $C_s(F_i), C_s(F_j)$.   The distance between two sub-cones, 
% $C_s(F_i) $ and $C_s(F_j) $ of PSD cone  $\mathbb{S}^n_+$, 
written as a function of frames $(F_i, F_j) $, is given by
\begin{align}\label{eq:Definition_Cone_Distance}
    d_{\mathrm{ cone}}(F_i,F_j) \!=\! \! \|F_iI_sF_i^T \!- \! F_jI_sF_j^T  \|_F.
\end{align}
\end{defn}
By definition of PSD matrices, for all sub-cones $F_k[ \mathbb{S}_+^s]F_k^T$, it always have the zero matrix inside it, $ 0 \in F_k[ \mathbb{S}_+^s]F_k^T$. Therefore, the restriction to identity matrix $ I_s$ is necessary to avoid trivial solutions, i.e., the zero matrix. 
The next theorem shows that the cone distance \eqref{eq:Definition_Cone_Distance} defined on PSD cone is equivalent to chordal distance \eqref{eq:Definition_Chordal_Distance} defined on the Grassmannian manifold $G(s, n)$, up to a scalar factor. 
% This will allow us then to leverage tools aim at (approximately) solving  the classical Grassmannian packing problem (Problem \ref{Prob:Grass_Packing}) to solve Problem \ref{Prob:PSD_Packing}.

\begin{thm} \label{Thm: Cone_Distance_equivalent}
The cone distance \eqref{eq:Definition_Cone_Distance} is equivalent to chordal distance \eqref{eq:Definition_Chordal_Distance} up to a scalar factor. More precisely, we have
\begin{align*}
 d_{\mathrm{ cone}}(F_i,F_j) = \sqrt{2}  d_{\mathrm{ chord}}(F_i,F_j) 
 \end{align*}
\textit{Proof}:  See Appendix~\ref{appendix:C}.

\end{thm}
% \textit{Remark:} Blekherman et al. introduced the \textit{normalized Frobenius distance} between the k-PSD closure $\mathbb{S}_+^{n,k} $and the PSD cone $\mathbb{S}_+^n$ : 
% \begin{align*}
%     \overline{\mathrm{dist}}_F(\mathbb{S}_+^{n,k},\mathbb{S}_+^n):= \sup_{X \in \mathbb{S}_+^{n,k}, \|X \|_F = 1 } \inf_{Y \in \mathbb{S}_+^n} \|X -Y \|_F
% \end{align*}
% We could apply this distance measure for the two sub-cones $ \overline{\mathrm{dist}}_F(C_s(F_i), C_s(F_j)):=$
% \begin{align*}
%     \sup_{X \in C_s(F_i), \|X \|_F = 1 } \inf_{Y \in C_s(F_j)} \|X -Y \|_F
% \end{align*}
% Applying first order necessary condition for optimality and some linear algebra, we have the following lemma establishing the connections between cone distance and normalized Frobenius distance
% \begin{lem} \label{lem:Frobneius Distance}
% $ \overline{\mathrm{dist}}_F(C_s(F_i), C_s(F_j))= [1 - \min_{k} \cos^4\theta_k ]^{1/2}$, where  $ \theta _k \in [0, \pi/2], k \in [s]$ denotes the principal angles  between  $ \mathcal{F}_i $ and $ \mathcal{F}_j  $.
% \end{lem}
% While the normalized Frobenius distance considers the $\min \cos^4\theta_k$, i.e., the largest principle angle $\theta_k$ between subspaces $ \mathcal{F}_i $ and $ \mathcal{F}_j  $, the cone distance takes all principal angles into account.

\begin{prob}[Subcone Packing in $\mathbb{S}^n_+$]\label{Prob:PSD_Packing}
Given the ambient PSD cone $\mathbb{S}^n_+$, and $s$-dimensional sub-cones of the form $F_k[ \mathbb{S}_+^s]F_k^T\}$. Find a set of N frames $\mathcal{F}= \{F_1, \dots,F_N \} \subseteq V_s(\eucd^n)$  that solves the mathematical program
\begin{align}
    \max_{\mathcal{F}: | \mathcal{F} |=N } \;\; \min_{F_i,F_j \in \mathcal{F}, i\neq j } d_{\mathrm{cone}}(F_i,F_j).
\end{align}
\end{prob}
At first sight, the question of finding frames that solves the mathematical program seems an unwieldy task, because of its highly nonconvex and combinatorial nature. However, as a direct result of Theorem \ref{Thm: Cone_Distance_equivalent}, the packing problem of  $s$-dimensional sub-cones in PSD cone in $\mathbb{S}^n_+$ (Problem \ref{Prob:PSD_Packing})  is equivalent to packing problem in Grassmannian manifold $G(s, n)$ (Problem \ref{Prob:Grass_Packing}). Therefore, we could leverage existing literature on constructing Grassmannian packings as a method of finding good sub-cone packings for approximating the PSD cone. When $N$ is fixed as in   $ \sum_{k=1}^N F_k[ \mathbb{S}_+^s]F_k^T $, these sub-cones are arranged so that they are as far apart possible and tend to provide more accurate inner approximations. By focusing on Grassmannian packings, the proposed approach provides a more flexible family of inner approximations as in \eqref{eq:Final_approximation}.

\section{Numerical Examples} \label{sec:numerical}
In this section, we leverage computational methods for finding good packings in Grassmannian manifolds via alternating projection, see  \cite{dhillon2008constructing}, to get build good approximations of the PSD cone. We illustrate the effectiveness of our proposed method with two numerical examples. The first one compares the inner approximation accuracy of scaled diagonally dominant ($SDD_n$) matrices, Factor-width 3 matrices ($\mathcal{FW}_n^3$) with our approach. The second numerical example studies how many sub-cones are needed to approximate the ambient PSD cone with high accuracy.
\subsection{Numerical Examples Setup}
We consider the following numerical methods for constructing Grassmannian packings, for its flexibility and ease of implementation. Consider a packing problem in Grassmannian Manifolds (Problem \ref{Prob:Grass_Packing}), suppose that the set of N $k-$dimensional subspaces $ \{\mathcal{S}_1,\dots,\mathcal{S}_N \} \subseteq G(k, n) $  are spanned by the matrices $\mathcal{F}=\{F_1,\dots,F_N \}$. Collate these $N$ matrices into a $n \times KN $ configuration matrix $F:= [F_1 \; \dots \;F_N ] $, and the block Gram matrix of $F$ is defined as the $  KN \times KN$ matrix $ G = F^*F$, whose blocks control the distances between pairs of subspaces. Solving the packing problem in Grassmannian manifold $G(k, n)$ is equivalent with constructing a Gram matrix $G$ that satisfies six structural and spectral properties, where the structural properties constrain the entries of the Gram matrix while the spectral properties control the eigenvalues of the matrix. The alternating projection algorithm alternately enforce the structural properties and then the spectral properties in hope of producing a Gram matrix that satisfies all six properties at once. Finally, factorization of the output Gram matrix $G = F^*F$, such as eigenvalue decomposition, yields the desired configuration matrix.
We note however some limitations of this approach, as pointed out in \cite{dhillon2008constructing}. First, although the alternating projection algorithm seems to converge in practice, a theoretical proof is lacking. Second, The algorithm frequently requires as many as 5000 iterations before the iterates settle down. Third, the algorithm was especially successful for smaller numbers of subspaces, but its performance began to flag as the number of subspaces approached 20.

In the following two numerical examples, the Matlab toolboxes YALMIP \cite{lofberg2004yalmip} and MOSEK \cite{mosek} are used for solving the semidefinite programming problems. 
\subsection{Approximation Error Comparison}
Consider the following semidefinite program to test the inner approximation accuracy of PSD cone:
\begin{align} \label{eq:Numerical}
    \min_{X \in \mathbb{S}^n}\;\; & \|X -A \|_F \nonumber \\
     s.t.\;\; &X \in {\mathcal{C}} \subseteq \mathbb{S}^n_+,
\end{align}
where $A$ is a random normalized PSD matrix, i.e., $\| A\|_F=1$ and $\mathcal{C}$ represents the sub-cones for different inner approximations. We calculate the distance between matrix $X$ and a random given PSD matrix $A$. The approximation quality is quantified as the empirical mean of optimal value, $\|X -A \|_F$, of \eqref{eq:Numerical}. 
% range from 0 to 1.

\begin{figure}[hbt!]
\centerline{\includegraphics[width=\columnwidth]{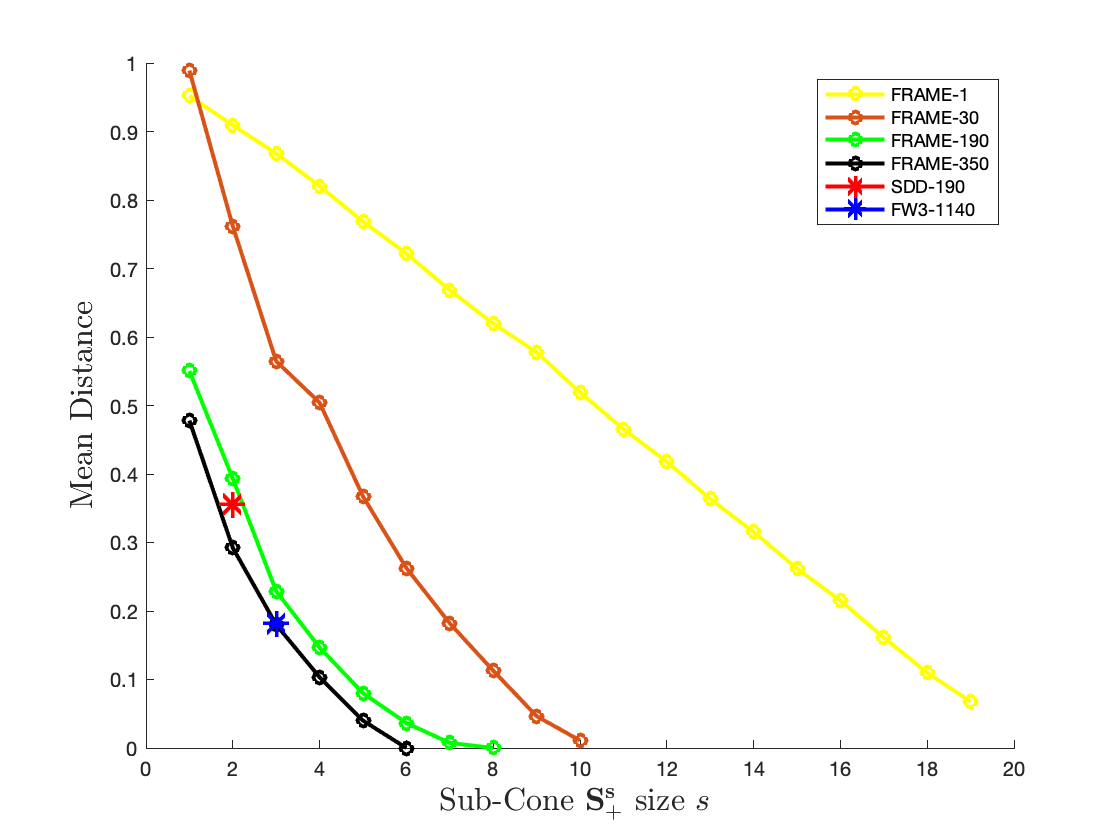}}
\caption{Comparison of the approximation error between scaled diagonally dominant (SDD), factor-width 3 matrices (FW3), and maximally distant frame set (FRAME) of N=$ \{1,30,190,350 \}$ with ambient cone $\mathbb{S}_+^{20}$ }
\label{fig:1 ApproxA}
\end{figure}

The simulation results are illustrated as in Figure \ref{fig:1 ApproxA}, with the ambient PSD cone of size $n = 20$. Here, the red curve $SDD-190$ denotes the case when  $\mathcal{C} = SDD_n$, and $190$ indicates the number of sub-cones $\binom{20}{2}=190$.
% comes from the fact that  any scaled diagonally dominant matrix $M$ can be written as $ M = \sum _{i=1}^{\binom{n}{2}} V_i \Lambda_i V_i^T $, i.e., enforcing the constraint $X \in SDD_n$ one needs $\binom{20}{2}=190$ sub-cones. 
Similarly, the dark blue curve $FW3-1140$ represents the approximations using $\binom{20}{3}=1140$ factor-width 3 matrices $\mathcal{C} = \mathcal{FW}_n^3$.
% and similarly one needs $\binom{20}{3}=1140$ sub-cones for this approach. 
The other three curves named by Frame-N denote a family of inner approximations  which takes the form $ \mathcal{C}=  \sum_{k=1}^N F_k[ \mathbb{S}_+^s]F_k^T$ as in \eqref{eq:Final_approximation}.
% , where $  F_k[\mathbb{S}^s_+ ]  F_k^T := \{ X \in \mathbb{S}^n_+ :  X =  F_k Y  F_k^T, Y \in \mathbb{S}^{s}_+  \}$. 
Here, four different packing have been found using the alternating projection algorithm~\cite{dhillon2008constructing}, where the numbers 1,30,190,350 represents the number N of sub-cones (equivalently number of frames). In each case, we gradually increase the sub-cone dimension $s$ to obtain tighter inner approximations of ambient PSD cone $\mathbb{S}^n_+$. 
The mean distance is calculated as the mean of 100 random tests. A more accurate inner approximation is then indicated by smaller mean distance. Notably, when the sub-cones are spanned by $\mathbb{S}_+^{10}$, one could reach almost perfect approximations with only $ N = 30$ frames, compared with $190$ sub-cones for SDD approach and $1140$ sub-cones for FW3 approach. Also, instead of using $1140$ sub-cones as in FW3 approach, we achieve a similar approximation accuracy with only $350$ sub-cones.

We end by noting that the SDD approach provides a more accurate inner approximation of the ambient PSD cone than the FRAME-190 approach when the sub-cone size $s=2$, i.e., when the two approaches both use 190 sub-cones of size $s=2$. We provide two possible explanations: 1) The $190$ frames $V_i$ as in  $ M = \sum _{i=1}^{\binom{20}{2}} V_i \Lambda_i V_i^T $ are (nearly) optimal packed; i.e., they are (nearly) maximally distant from each other and thus provides the (nearly) optimal approximation; 2) As mentioned in above subsection, the performance of the alternating projection algorithm used in FRAME approach begins to flag as the number of subspaces approached 20. Since the number of sub-cones is 190, the algorithm could not construct optimal packings and cannot provide the most accurate approximations.

\subsection{Number of Sub-cones}
\begin{figure}[hbt!]
\centerline{\includegraphics[width=\columnwidth]{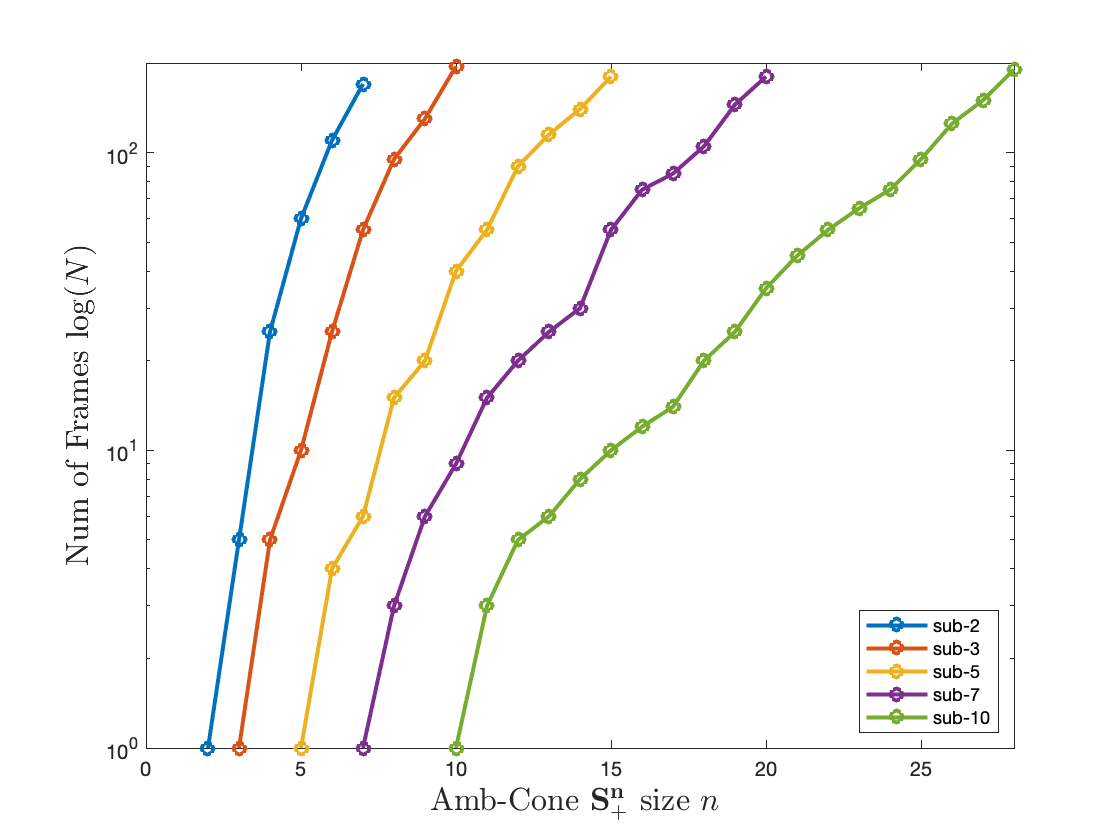}}
\caption{Calculating number of frames needed for a tight inner approximation of the ambient PSD cone. Here, the number $s$ in sub-s represents the dimension of sub-cone $\mathbb{S}_+^{s} $ }
\label{fig:2 ApproxB}
\end{figure}

In the second numerical example, we explore the number of sub-cones (frames) required for a tight inner approximation of the ambient PSD cone. Specifically, we consider the same SDP problem as in \eqref{eq:Numerical}
\begin{align*}
    \min_{X \in \mathbb{S}^n}\;\; & \|X -A \|_F \nonumber \\
     s.t.\;\; &X \in {\mathcal{C}} \subseteq \mathbb{S}^n_+,
\end{align*}
where $ \mathcal{C}=  \sum_{k=1}^N F_k[ \mathbb{S}_+^s]F_k^T$; i,e., optimization variable $X$ is approximated as sum of sub-cones. 
The simulation results are illustrated as in Figure \eqref{fig:2 ApproxB}. Here, we are looking for the numbers of sub-cones needed in each case such that the average mean distance $\|X -A \|_F$ through 100 random tests is less than $0.01$. The $x$ axis represents the dimension of ambient PSD cone $n$ and the y axis represents the numbers of 
frames required in logarithmic scale. For example, the blue curve sub-2 represents the case when the sub-cones are spanned by $ \mathbb{S}_+^{2}$. Therefore, it is always possible to cover the ambient PSD cone $\mathbb{S}_+^{2}$ with only one sub-cone of the same size $\mathbb{S}_+^{2}$. From the simulation, the curves appears linearly in the logarithmic scale, which indicates the number of frames required for a tight inner approximation grows exponentially.  Note that this simulation result is consistent with the result in \cite{song2021approximations}, which proved that the number of subspaces must be at least exponentially large in $n$ for a good approximation of the PSD cone $\mathbb{S}^n_+$.

\section{Conclusions}\label{sec:conlusions}
In this paper, we developed a novel decomposition framework of the PSD cone by means of conic combinations of sub-cones, which unified many existing inner approximation techniques. Furthermore, we introduced a more flexible family of inner approximations of the PSD cone by a set of sub-cones, where we aimed to arrange these sub-cones so that they are maximally separated from each other. We showed the problem of packing sub-cones is equivalent to a packing problem in Grassmannian manifolds with chordal distances. The effectiveness of our approach was demonstrated with simulation results.

\appendix

\subsection{Proof of Theorem \ref{lem:General PSD Decomposition Framework}}\label{appendix:A}

Given any PSD matrix $ X \in \mathbb{S}^n_+$, it is diagonalizable by some $ U \in O(n)$ , i.e., $ X = U \Lambda U^T = \sum_{i=1}^n \lambda_i u_i u_i^T$, where $u_i$ is the $i$th column of $ U$ and $\lambda_i $ is the $i$th diagonal element of $ \Lambda $. Let $ S = \{ 1 \}$, and $ \mathcal{F} \subseteq O(n) :=\{ F_k, k \in  [n]\}$ such that the first column of $ F_k$ equal to $u_k$, i.e., the $k$th column of $U$. By construction,  we have $X \in \mathrm{cone} (\bigcup_{ F \in \mathcal{F}} F_S[\mathbb{S}^{1 }_+ ]  F_S^T  )$.

Conversely, consider any element $Y_i \in \mathbb{S}^n$ of the set $ \bigcup_{S \in \mathcal{S} } \bigcup_{ F \in O(n)} F_S[\mathbb{S}^{|S| }_+ ]  F_S^T $, i.e., $Y_i =F_S Y  F_S^T  $ for some $S \in \mathcal{S} $, $ F \in O(n)$ and $ Y \in \mathbb{S}^{|S| }_+ $. Firstly, we have $ x^T Y_i x = x^T F_S Y  F_S^T x \geq 0 ,\; \forall x \in \eucd ^n$, since $Y$ is a $ |S| \times |S| $ PSD matrix and $F_S^T x \in \eucd^{|S|}$. By definition, (possibly infinite) sum of PSD matrices are also PSD matrices. Next, we want to show the conical combination $ Z = \sum _{i=1}^k \alpha_i Y_i $, for some $ Y_i \in  \bigcup_{S \in \mathcal{S} } \bigcup_{ F \in O(n)} F_S[\mathbb{S}^{|S| }_+ ]  F_S^T $, with $ \alpha_i \geq 0, i \in [k], k \in \mathbb{N}$, being PSD. Again, applying the quadratic form $ x^T Z x = \sum _{i=1}^k \alpha_ix^T Y_i x \geq 0  $, which completes the proof. 

\subsection{Proof of Theorem \ref{lem:simplified_Decomposition}} \label{appendix:B}

The proof follows the same structure as in Theorem \ref{lem:General PSD Decomposition Framework}.
Given any PSD matrix $ X \in \mathbb{S}^n_+$, its diagonalization is given by $ X = U \Lambda U^T = \sum_{i=1}^n \lambda_i u_i u_i^T$. Let $ X_i \in \mathbb{S}^s _+ , i\in [n]$ be a diagonal matrix with $ (X_i)_{11} = \lambda_i\|u_i \|_2^2$ and other being zero. $F_i \in V_s(\eucd^n ), i\in [n]$ such that the first column equals $u_i/ \|u_i \|_2$. Thus we have $ X = \sum_{i=1}^n F_i X_i F_i^T$.

Conversely, consider any element $Y_i \in \mathbb{S}^n$ of the set $ \bigcup_{ F \in V_s(\eucd^n )} F[\mathbb{S}^s_+ ]  F^T $ for some $ F \in  V_s(\eucd^n )$ and $ Y \in \mathbb{S}^{s }_+ $. Firstly, we have $ x^T Y_i x\geq 0 ,\; \forall x \in \eucd ^n$. By definition, (possibly infinite) sum of PSD matrices are also PSD matrices. Next, we want to show the conical combination $ Z = \sum _{i=1}^k \alpha_i Y_i $, for some $ Y_i \in  \bigcup_{ F \in V_s(\eucd^n )} F[\mathbb{S}^s_+ ]  F^T$, with $ \alpha_i \geq 0, i \in [k], k \in \mathbb{N}$, being PSD. Again, applying the quadratic form $ x^T Z x = \sum _{i=1}^k \alpha_ix^T Y_i x \geq 0  $, which completes the proof.

\subsection{Proof of Theorem \ref{Thm: Cone_Distance_equivalent}} \label{appendix:C}
Since $Y_i , Y_j   \in   O(s)$, we could simplify the cone distance  \eqref{eq:Definition_Cone_Distance} as follows
\begin{align*}
d_{\mathrm{ cone}}(F_i,F_j)^2 &= \| F_iF_i^T -F_jF_j^T \|_F^2 \\
&= \mathrm{trace}(( F_iF_i^T -F_jF_j^T )^T (F_iF_i^T -F_jF_j^T )) \\
&= \mathrm{ trace}( F_iF_i^T F_iF_i^T )+ \mathrm{ trace}( F_jF_j^T F_jF_j^T ) \\
 &-2 \mathrm{ trace}(F_iF_i^T F_jF_j^T) \\
&= 2s - 2 \| F_i^T F_j \|_F^2 \\
&= 2 d_{\mathrm{ chord}}(F_i,F_j)^2
\end{align*}
where we used the fact that $F_i^T F_i = I_s $ and the cyclic property of trace operator.

\bibliographystyle{IEEEtran}

\bibliography{main.bib} 
\end{document}

%% file: edits.tex
\usepackage{ifthen}
\newboolean{showcomments}
\setboolean{showcomments}{true}
\usepackage[usenames,dvipsnames]{xcolor}
\usepackage{todonotes}

\definecolor{bleudefrance}{rgb}{0.19, 0.55, 0.91}
\definecolor{ao(english)}{rgb}{0.0, 0.5, 0.0}

\newcommand{\addcite}[0]{\ifthenelse{\boolean{showcomments}}
{\textcolor{purple}{(add cite(s)) }}{}}%

\newcommand{\enrique}[1]{  \ifthenelse{\boolean{showcomments}}
{\todo[inline,color=bleudefrance]{Enrique: #1}}{}}
\newcommand{\emmargin}[1]{\ifthenelse{\boolean{showcomments}}{\marginpar{\color{bleudefrance}\tiny EM: #1}}{}}

\newboolean{showedits}
\setboolean{showedits}{true}
\usepackage[xcolor=bleudefrance]{changes}
\definechangesauthor[color=bleudefrance]{EM}
\newcommand{\aem}[1]{
\ifthenelse{\boolean{showedits}}
{\added[id=EM]{#1}}
{#1}
}
\newcommand{\chem}[2]{
\ifthenelse{\boolean{showedits}}
{\replaced[id=EM]{#1}{#2}}
{#1}
}
\newcommand{\dem}[1]{
\ifthenelse{\boolean{showedits}}
{\deleted[id=EM]{#1}}
{}
}

%% file: main.bbl
\begin{thebibliography}{10}
\providecommand{\url}[1]{#1}
\csname url@rmstyle\endcsname
\providecommand{\newblock}{\relax}
\providecommand{\bibinfo}[2]{#2}
\providecommand\BIBentrySTDinterwordspacing{\spaceskip=0pt\relax}
\providecommand\BIBentryALTinterwordstretchfactor{4}
\providecommand\BIBentryALTinterwordspacing{\spaceskip=\fontdimen2\font plus
\BIBentryALTinterwordstretchfactor\fontdimen3\font minus
  \fontdimen4\font\relax}
\providecommand\BIBforeignlanguage[2]{{%
\expandafter\ifx\csname l@#1\endcsname\relax
\typeout{** WARNING: IEEEtran.bst: No hyphenation pattern has been}%
\typeout{** loaded for the language `#1'. Using the pattern for}%
\typeout{** the default language instead.}%
\else
\language=\csname l@#1\endcsname
\fi
#2}}

\bibitem{blekherman2012semidefinite}
G.~Blekherman, P.~A. Parrilo, and R.~R. Thomas, \emph{Semidefinite optimization
  and convex algebraic geometry}.\hskip 1em plus 0.5em minus 0.4em\relax SIAM,
  2012.

\bibitem{vandenberghe1996semidefinite}
L.~Vandenberghe and S.~Boyd, ``Semidefinite programming,'' \emph{SIAM review},
  vol.~38, no.~1, pp. 49--95, 1996.

\bibitem{parrilo2000structured}
P.~A. Parrilo, ``Structured semidefinite programs and semialgebraic geometry
  methods in robustness and optimization,'' Ph.D. dissertation, California
  Institute of Technology, 2000.

\bibitem{8431617}
L.~{Wang}, D.~{Han}, and M.~{Egerstedt}, ``Permissive barrier certificates for
  safe stabilization using sum-of-squares,'' in \emph{2018 Annual American
  Control Conference (ACC)}, 2018, pp. 585--590.

\bibitem{deits2015computing}
R.~Deits and R.~Tedrake, ``Computing large convex regions of obstacle-free
  space through semidefinite programming,'' in \emph{Algorithmic foundations of
  robotics XI}.\hskip 1em plus 0.5em minus 0.4em\relax Springer, 2015, pp.
  109--124.

\bibitem{5717711}
J.~{Derenick}, J.~{Spletzer}, and V.~{Kumar}, ``A semidefinite programming
  framework for controlling multi-robot systems in dynamic environments,'' in
  \emph{49th IEEE Conference on Decision and Control (CDC)}, 2010, pp.
  7172--7177.

\bibitem{benson2000solving}
S.~J. Benson, Y.~Ye, and X.~Zhang, ``Solving large-scale sparse semidefinite
  programs for combinatorial optimization,'' \emph{SIAM Journal on
  Optimization}, vol.~10, no.~2, pp. 443--461, 2000.

\bibitem{majumdar2019survey}
A.~Majumdar, G.~Hall, and A.~A. Ahmadi, ``A survey of recent scalability
  improvements for semidefinite programming with applications in machine
  learning, control, and robotics,'' \emph{arXiv preprint arXiv:1908.05209},
  2019.

\bibitem{ahmadi2019dsos}
A.~A. Ahmadi and A.~Majumdar, ``Dsos and sdsos optimization: more tractable
  alternatives to sum of squares and semidefinite optimization,'' \emph{SIAM
  Journal on Applied Algebra and Geometry}, vol.~3, no.~2, pp. 193--230, 2019.

\bibitem{zheng2019block}
Y.~Zheng, A.~Sootla, and A.~Papachristodoulou, ``Block factor-width-two
  matrices and their applications to semidefinite and sum-of-squares
  optimization,'' \emph{arXiv preprint arXiv:1909.11076}, 2019.

\bibitem{blekherman2020sparse}
G.~Blekherman, S.~S. Dey, M.~Molinaro, and S.~Sun, ``Sparse psd approximation
  of the psd cone,'' \emph{Mathematical Programming}, pp. 1--24, 2020.

\bibitem{song2021approximations}
D.~Song and P.~A. Parrilo, ``On approximations of the psd cone by a polynomial
  number of smaller-sized psd cones,'' \emph{arXiv preprint arXiv:2105.02080},
  2021.

\bibitem{boman2005factor}
E.~G. Boman, D.~Chen, O.~Parekh, and S.~Toledo, ``On factor width and symmetric
  h-matrices,'' \emph{Linear algebra and its applications}, vol. 405, pp.
  239--248, 2005.

\bibitem{vandenberghe2015chordal}
L.~Vandenberghe and M.~S. Andersen, ``Chordal graphs and semidefinite
  optimization,'' \emph{Foundations and Trends in Optimization}, vol.~1, no.~4,
  pp. 241--433, 2015.

\bibitem{conway1996packing}
J.~H. Conway, R.~H. Hardin, and N.~J. Sloane, ``Packing lines, planes, etc.:
  Packings in grassmannian spaces,'' \emph{Experimental mathematics}, vol.~5,
  no.~2, pp. 139--159, 1996.

\bibitem{kutyniok2009robust}
G.~Kutyniok, A.~Pezeshki, R.~Calderbank, and T.~Liu, ``Robust dimension
  reduction, fusion frames, and grassmannian packings,'' \emph{Applied and
  Computational Harmonic Analysis}, vol.~26, no.~1, pp. 64--76, 2009.

\bibitem{love2003grassmannian}
D.~J. Love, R.~W. Heath, and T.~Strohmer, ``Grassmannian beamforming for
  multiple-input multiple-output wireless systems,'' \emph{IEEE transactions on
  information theory}, vol.~49, no.~10, pp. 2735--2747, 2003.

\bibitem{gohary2009noncoherent}
R.~H. Gohary and T.~N. Davidson, ``Noncoherent mimo communication: Grassmannian
  constellations and efficient detection,'' \emph{IEEE Transactions on
  Information Theory}, vol.~55, no.~3, pp. 1176--1205, 2009.

\bibitem{strohmer2003grassmannian}
T.~Strohmer and R.~W. Heath~Jr, ``Grassmannian frames with applications to
  coding and communication,'' \emph{Applied and computational harmonic
  analysis}, vol.~14, no.~3, pp. 257--275, 2003.

\bibitem{dhillon2008constructing}
I.~S. Dhillon, J.~R. Heath, T.~Strohmer, and J.~A. Tropp, ``Constructing
  packings in grassmannian manifolds via alternating projection,''
  \emph{Experimental mathematics}, vol.~17, no.~1, pp. 9--35, 2008.

\bibitem{bjorck1973numerical}
A.~Björck and G.~H. Golub, ``Numerical methods for computing angles between
  linear subspaces,'' \emph{Mathematics of computation}, vol.~27, no. 123, pp.
  579--594, 1973.

\bibitem{alizadeh1995interior}
F.~Alizadeh, ``Interior point methods in semidefinite programming with
  applications to combinatorial optimization,'' \emph{SIAM journal on
  Optimization}, vol.~5, no.~1, pp. 13--51, 1995.

\bibitem{barker1975cones}
G.~Barker and D.~Carlson, ``Cones of diagonally dominant matrices,''
  \emph{Pacific Journal of Mathematics}, vol.~57, no.~1, pp. 15--32, 1975.

\bibitem{lofberg2004yalmip}
J.~Lofberg, ``Yalmip: A toolbox for modeling and optimization in matlab,'' in
  \emph{2004 IEEE international conference on robotics and automation (IEEE
  Cat. No. 04CH37508)}.\hskip 1em plus 0.5em minus 0.4em\relax IEEE, 2004, pp.
  284--289.

\bibitem{mosek}
\BIBentryALTinterwordspacing
M.~ApS, \emph{The MOSEK optimization toolbox for MATLAB manual. Version 9.0.},
  2019. [Online]. Available: \url{http://docs.mosek.com/9.0/toolbox/index.html}
\BIBentrySTDinterwordspacing

\end{thebibliography}
